\documentclass{amsart}

\usepackage{graphicx}

\def\R{{\mathbb R}}

\def\Z{{\mathbb Z}}
\def\1{{1\!\!\!1}}

\def\E{{\mathbb E}}
\def\P{{\mathbb P}}

\def\cal{\mathcal}

\def\eps{\varepsilon}

\newcommand{\be}{\begin{equation}}
\newcommand{\ee}{\end{equation}}
\numberwithin{equation}{section}

\newtheorem{theorem}{Theorem}
\newtheorem{prop}{Proposition}[section]
\newtheorem{cor}{Corollary}[section]

\newtheorem{lemma}{Lemma}[section]

\title{Quasipotential and  logarithmic asymptotics of  Green's measures}
\author{Irina Ignatiouk-Robert}
\address{
{Universit\'e de Cergy-Pontoise,}
{D\'epartement de math\'ematiques,}
{2, Avenue Adolphe Chauvin,}
{95302 Cergy-Pontoise Cedex,}
{France}}
\date{\today}
\email{Irina.Ignatiouk@math.u-cergy.fr}

\begin{document}
\begin{abstract} It is proved 
that the weak large deviation principle of the scaled processes $Z^\eps(t)=\eps Z(t/\eps)$
implies the weak large deviation principle for the scaled Green's measures 
of the  Markov process $Z(t)$.  
\end{abstract}
\maketitle

\section{Introduction and main results}\label{sec1}

For a random perturbation
of a dynamical system (see Freidlin and Wentzell~\cite{W-F}) satisfying sample path large
deviation principle, quasipotential  of the
corresponding rate function $I_{[T_1,T_2]}$ on the Skorohod space $D([T_1,T_2],\R^d)$, $T_1\leq T_2$,
  is defined by the equality 
\[
I(q,q') ~\dot=~ \inf_{0 \leq T_1 < T_2< \infty} ~\inf_{\phi : \phi(T_1)=q, \phi(T_2)=q'}
~I_{[T_1,T_2]}(\phi). 
\]
This function is important in several problems. Quasipotential  characterises the asymptotical
behavior of the stationary probabilities and also of the expectation of first exit time
from the domain. Intuitively, the quantity $I(q,q')$ represents an optimal large deviation
cost of going from the small  
neighborhood of $q$ to a small neighborhood of $q'$ within a reasonable time.

For scaled homogeneous random walks  $S^\eps(t) = \eps S([t/\eps])$, the quasipotential is given by  
\[
I(0,q)  ~=~ \sup_{a\in\R^d~: \varphi(a)\leq 1} ~a\cdot q 
\] 
where $\varphi$ is a jump generating function defined by  
\[
\varphi(a) ~\dot=~ \E_0\left(e^{a\cdot S(1)}\right). 
\]
For these processes, the function $I(0,q)$ characterizes the
asymptotical behavior of the Green's (potential) function  
\[
G(0,y) ~\dot=~ \sum_{t=0}^\infty \P_0(S(t) = y).
\] 
From the exact asymptotics obtained by  Ney and
Spitzer~\cite{Ney-Spitzer} and Borovkov and
Mogulskii~\cite{Borovkov-Mogulskii:02}  it follows a weak large
deviation principle with the rate function $I(0,q)$ :  
\[
\liminf_{n\to\infty} \frac{1}{n} \log G(0, nO) ~\geq~ - \inf_{q\in O} I(0,q) 
\]
for every open set $O\subset\R^d$ and 
\[
\limsup_{n\to\infty} \frac{1}{n} \log G(0, nV) ~\geq~ - \inf_{q\in V} I(0,q) 
\]
for every compact set $V\subset\R^d$. Remark that the hole large deviation principle does
not hold because $G(0,\R^d) = \infty$.

\medskip

In the present paper, such a property is extended for non homogeneous  Markov processes 
$Z(t)$ on $\R^d$: it is proved that  if the Markov process $Z(t)$ is transient and the sequence of
scaled processes $Z_n(t) = 
Z(nt)/n$ on the interval $[0,T]$ satisfies
sample path large deviation principle with a good rate function $I_{[0,T]}$ and if 
\[
I_T(0,0) ~\dot=~ \inf_{\phi : \phi(0)=\phi(T)=0}
~I_{[0,T]}(\phi) ~>~ 0 
\]
then for any $q_0\in\R^d$ and any sequence of points $z_n\in \R^d$ with $\lim_{n\to\infty}
z_n/n =q_0$, the
sequence of measures  
\[
\mu_n(B) ~=~ G(z_n, n B) ~\dot=~ \sum_{t=0}^\infty \P_{z_n}(Z(t)\in n B) 
\]
satisfies weak large deviation principle with
the rate function 
\[
I(q_0, q) ~=~ \inf_{T>0} ~\inf_{\phi : \phi(0)=q_0, \phi(T)=q} ~I_{[0,T]}(\phi).  
\]

This result is motivated by applications to the problem
of Martin boundary for partially homogeneous random walks for which sample path large
deviation principle for the sequence of scaled processes was obtained and the
corresponding rate function was identified while the Martin boundary was described
only in very particular cases (see the papers of Alili and Doney~\cite{Alili-Doney},
Kurkova and Malyshev~\cite{Kurkova-Malyshev} and Ignatiouk~\cite{Ignatiouk:06}). In
general, in order to
identify the Martin compactification, one should calculate the exact asymptotics of the
Green's function. The  weak large deviation
principle provides the rough logarithmic asymptotics for the Green's function which  is
the first step in this direction.  Moreover, in some cases (see
Ignatiouk~\cite{Ignatiouk:06}), the rough logarithmic
asymptotics of the Green's function allow to describe the Martin
boundary in a straightforward way.

\subsection{Main result}
We consider  a strong  Markov process $(Z(t))$ on $E\subset\R^d$ whose sample
paths are right continuous and have the left limits.   To simplify the notations, it
is convenient to consider  continuous time Markov chains. For discrete time Markov chains,
all our results can be extended in a straightforward way, by replacing the variables
$Z(t)$ for $t\in\R_+$ by $Z([t])$ 
where $[t]$ denotes the integer part of $t$. 

The set $E$ is assumed to be
unbounded and the Green's function 
\[
G(z, B) ~\dot=~ \int_0^\infty \P_z(Z(t)\in B)  ~\, dt 
\]
is assumed to be well defined and finite for every $z\in E$ and every compact set
$B\subset \R^d$. For the Markov process $(Z(t))$ we consider a family of scaled processes
$Z^\eps(t)$ generated by $(Z(t))$ : 
\[
(Z^\eps(t), \; t\in[0,T]) ~\dot=~ (\eps Z(t/\eps),
\; t\in[0,T])
\]  
The Markov process $(Z(t))$ will be assumed  to satisfy the
following conditions :  

\medskip 

\noindent
{\bf (H1)} {\em Large deviations.}   For every $T>0$, the family of
rescaled processes $(Z^\eps(T))$ satisfies weak large deviation principle in $\R^d$ with
a good rate functions $I_T :\R^d\times\R^d\to\R_+$ :
\begin{itemize}
\item[--] the function $I_T :\R^d\times\R^d\to\R_+$ is lower semicontinuous;
\item[--] for any $q\in\R^d$ and any open set $O\subset\R^d$ 
\be\label{e-lb}
\lim_{\delta\to 0}~\liminf_{\eps\to 0} ~\eps \log ~\inf_{z\in E : |\eps z - q| < \delta}
~\P_z(Z^\eps(T) \in O) ~\geq~ - \inf_{q'\in O} I_T(q,q') 
\ee
\item[--] for any $q\in\R^d$ and any compact set $V\subset\R^d$ 
\be\label{e-ub}
\lim_{\delta\to 0}~\liminf_{\eps\to 0} ~\eps \log ~\sup_{z\in E : |\eps z - q| < \delta}
~\P_z(Z^\eps(T) \in V) ~\leq~ - \inf_{q'\in V} I_T(q,q') 
\ee 
\end{itemize} 
Here and throughout we denote by $\P_{z}$  a conditional probability given that $Z(0)=
z$.

\medskip 

\noindent 
{\bf (H2)} {\em Asymptotically finite range.} The function 
\[
\hat{\varphi}(a) ~\dot=~ ~\sup_{z\in E} ~\sup_{t\in[0,1]}  \E_z\left(e^{a\cdot
  (Z(t)-z)}\right)   
\] 
is finite everywhere on $\R^d$. 

\medskip

\noindent 
{\bf (H3)} {\em  Communication condition.} There are $\theta>0$ and positive function
  $\sigma : E\to\R_+$ such that $\sigma(z)/|z| \to
  0$ when $|z|\to \infty$ and for
  every $z,z'\in E$, the 
  probability that starting at $z$, the 
  Markov chain $Z(t)$ ever hits the open ball $B(z',\sigma(z'))$ centered at $z'$ and
  having the radius $\sigma(z')$ is greater than $\exp(-\theta|z'-z|)$.

\medskip 

Remark that by contraction principle, the condition (H1) is satisfied under  the following
assumption.  

\medskip 

\noindent
{\bf (H1')} {\em Sample path large deviations.}   For every $T>0$, the family of
scaled processes $(Z^\eps(t), \; t\in[0,T])$ satisfies sample 
path large deviation principle in  the Skorohod space $D([0,T],\R^d)$ with a good rate functions
$I_{[0,T]}$ and  for any $q,q'\in\R^d$,
$I_T(q,q')$ is the infimum
of the rate function $I_{[0,T]}(\phi)$ over all $\phi\in D([0,T],\R^d)$ with $\phi(0)=q$
and $\phi(T)=q'$ (see section~\ref{sec2} for more details). 

\medskip

Let ${\cal R}$ denote the set of all possible limits $\lim_{\eps\to 0} \eps z_\eps$ with
$z_\eps\in E$. For given $q,q'\in{\cal R}$  we let 
\[
I(q,q') ~\dot=~ \inf_{T>0} ~I_T(q,q') \quad \text{ and } \quad \hat{I}(q,q') ~\dot=~
\begin{cases} I(q,q') &\text{if $q\not= q$}\\ 0 &\text{if $q=q'$} \end{cases}
\]
It is convenient moreover to introduce the following notations : for $R>0$ we
let 
\[
\tau_R ~\dot=~ \inf\{ t\geq 0 :~ |Z(t)| \geq R\} 
\]
and we consider the truncated Green's function 
\[
G_R(z, B) ~=~ \int_0^\infty \P_z\left( Z(t) \in B, \; \tau_R > t\right) \, dt.
\]

\bigskip 
\noindent 
Our main result is the following theorem. 

\begin{theorem}\label{th1} Suppose that the conditions (H1) and (H2) are satisfied and let
 $I_T(0,0) > 0$. Then 
the following assertions hold : 
\begin{itemize}
\item[(i)] for any $q\in{\cal R}$ and any open set $O\subset\R^d$, 
\be\label{e1-1}
\lim_{\delta\to 0} ~\liminf_{n\to\infty} ~\frac{1}{n}\log  ~\inf_{z\in E : |z- n q| < 
  \delta a} G\left(z, n O\right) ~\geq~ - \inf_{q'\in O} I(q,q'), 
\ee
\item[(ii)] for any $q\in{\cal R}$ and any  compact set $V\subset\R^d$, 
\be\label{e1-2}
\lim_{\delta\to 0} ~\limsup_{n\to\infty} ~\frac{1}{n}\log  ~\sup_{z\in E : |z-n q| < 
  \delta a} G(z,  n V) ~\leq~ - \inf_{q'\in V} \hat{I}(q,q') 
\ee
\end{itemize}
If the assumption (H1') is also satisfied then 
\begin{itemize}
\item[(iii)] for any  bounded set $V\subset\R^d$, any $q\in{\cal R}$ and for any
  $A>0$ there is $R>0$ such that  
\be\label{e1-3}
\lim_{\delta\to 0} ~\limsup_{n\to \infty} ~\frac{1}{n}\log  ~\sup_{z\in E : |z-n q| < 
  \delta a} \bigl(G(z, n V) - G_{nR}(z, nV)\bigr) ~\leq~  - A 
\ee
\end{itemize} 
If the conditions (H1) - (H3) are satisfied with $I_T(0,0) >0$ then 
\begin{itemize}
\item[(iv)] $I(q,q')\equiv \hat{I}(q,q')$ and
the function $I(q,q')$ is continuous on ${\cal R}\times{\cal R}$.  
\end{itemize} 
\end{theorem}
Remark that the first assertion of this theorem implies the  large deviation lower bound
for the sequence of measures $\mu_n(B) = G(z_n, nB)$ in its usual form : 
for any $q\in{\cal R}$ and any open set $O\subset\R^d$, 
\[
\lim_{\delta\to 0} ~\liminf_{n\to \infty} ~\frac{1}{n}\log  G(z_n, n O) ~\geq~ - \inf_{q'\in O} I(q,q')
\]
when $z_n/n \to q$ as $n\to\infty$. The second assertion of this theorem implies the 
upper bound on compact sets : for any $q\in{\cal R}$ and any compact set $V\subset\R^d$, 
\[
\lim_{\delta\to 0} ~\limsup_{n\to \infty} ~\frac{1}{n}\log  G(z_n, n V) ~\leq~ - \inf_{q'\in V} I(q,q')
\]
when $z_n/n \to q$ as $n\to\infty$. The third assertion shows that the main
contribution to the quantity $G(z_n, nV)$ is given by the probability of those trajectories of the
process $(Z(t))$ which do not exit from the open ball $B(0, n R)$ centered at $0$ and
having the radius $nR$. 

\subsection{Outline of the proof} The proof of the lower bound \eqref{e1-1} is
straightforward. For discrete time Markov process $(Z(t))$, the lower bound \eqref{e1-1}
follows from the lower large
deviation bound \eqref{e-lb} applied with the large deviation parameter $\eps=1/a$ and the inequality 
\[
G\bigl(z, n O \bigr) ~\geq~ \P_z\left(Z([Tn]) \in O \right) ~=~
\P_z\left(Z^{1/n}(T) \in O\right).    
\]
For continuous time Markov process $(Z(t))$,  we show that  for any $q'\in O$ and $\delta'>0$
such that $B(q',\delta')\subset O$ the following inequality holds   
\begin {align*}
\liminf_{n\to \infty} ~\frac{1}{n}\log  &~\inf_{z\in E : |z-nq| < 
  \delta n}  G\bigl(z, n O \bigr) \\ &~\geq~ \liminf_{n\to \infty}
  ~\frac{1}{n}\log  ~\inf_{z\in E : |z-nq| <  \delta n}   \P_z\left(Z^{1/n}(T) \in
  B(q',\delta'/2)\right).  
\end{align*} 
The lower large deviation bound \eqref{e-lb} is applied then with an open set
$B(q',\delta')$ for every $q'\in O$.

The proof of the upper bound \eqref{e1-2} is more technical. We show that for any $q,q'\in{\cal R}$ 
\be\label{e1-6}
\lim_{\delta\to 0} ~\limsup_{n\to \infty} ~\frac{1}{n}\log  ~\sup_{z\in E : |z-nq| < 
  \delta a} G(z, n B(q',\delta)) ~\leq~ -  \hat{I}(q,q').  
\ee 
Here, the straightforward application of the upper large
deviation bound \eqref{e-ub} would imply that for any $T>0$ and $\sigma >0$, there are $n_T > 0$ and
$\delta_T > 0$ such that  for $t=Tn$, 
\[
~\frac{1}{n}\log  ~\sup_{z\in E : |z-nq| < 
  \delta n} \P_z\bigl( Z(t) \in n B(q',\delta)\bigr)  ~\leq~ -  I_T(q,q') +  \sigma ~\leq~ -
\hat{I}(q,q') + \sigma 
\]
for all $n \geq n_T$ and $0< \delta < \delta_T$. These estimates are not sufficient for the
proof of \eqref{e1-6} because the number $n_T$ depends on $T$ (remark that  in our setting, the
function $T\to n_T$ is implicit). To get \eqref{e1-6} we change the scale : 
the upper large deviation bound \eqref{e-ub} is now used with the large deviation parameter
$\eps=1/t$. For $\kappa >0$ small enough, the upper bound 
\be\label{e1-7}
\lim_{\delta\to 0} ~\limsup_{n\to \infty} ~\frac{1}{n}\log  ~\sup_{z\in E : |z-nq| < 
  \delta n}  \int_0^{\kappa n} \P_z\bigl( Z(t) \in n B(q',\delta)\bigr) \, dt ~\leq~ -
  \hat{I}(q,q') 
\ee 
is proved by using Chebyshev's inequality. The upper large deviation bound 
\be\label{e1-8}
\lim_{\delta\to 0} ~\limsup_{\eps\to 0} ~\sup_{z\in E : |\eps z - x| < 
  \delta } ~\eps \log   ~\P_z\bigl(Z^{\eps}(1) \in B(y,\delta)\bigr) ~\leq~ -  I_1(x,y)   
\ee
with $x=y=0$ is used for the
proof of the inequality   
\be\label{e1-9}
\lim_{\delta\to 0} ~\limsup_{n\to \infty} ~\frac{1}{n}\log  ~\sup_{z\in E : |z-nq| < 
  \delta n}  \int_{K n}^\infty \P_z\bigl( Z(t) \in n B(q',\delta)\bigr) \, dt ~\leq~ -
  \hat{I}(q,q')  
\ee
for $K > 0$ large enough.  To prove the upper bound 
\be\label{e1-10}
\lim_{\delta\to 0} ~\limsup_{n\to \infty} ~\frac{1}{n}\log  ~\sup_{z\in E : |z-nq| < 
  \delta n}  \int_{\kappa n}^{K n} \P_z\bigl( Z(t) \in n B(q',\delta)\bigr) \, dt ~\leq~ -
  \hat{I}(q,q')  
\ee  
we use 
the inequality 
\[
~\P_z\Bigl( Z(t) \in n B(q',\delta)\Bigr) ~\leq~ \P_z\Bigl( Z^{1/t}(1) \in  B(n q'/t,\delta/\kappa)\Bigr) 
\]
and the upper large deviation bound \eqref{e1-8} with $\eps = 1/t$, $x= \theta q$
and $y=\theta q'$ for 
each  $\theta = n/t \in [K^{-1},\kappa^{-1}]$. 

The proof of the third assertion uses the inequality \eqref{e1-9} and the sample path
large deviation upper bound. To prove the last assertion of Theorem~\ref{th1} we combine
the upper bound \eqref{e1-2} and the rough lower bound  
\[
\lim_{\delta\to 0} ~\liminf_{n\to \infty} ~\frac{1}{n}\log  ~\inf_{z\in E : |z-nq| < 
  \delta a} G(z, n B(q',\delta))  ~\geq~ - \theta |q'-q| 
\]
which is a consequence  of the communication condition (H3). 

Our paper is organized as follows. In Section~\ref{sec2} we recall the large deviation
properties of the scaled processes. The inequalities \eqref{e1-7}, \eqref{e1-9} and
\eqref{e1-10} are proved in Section~\ref{sec3}. The proof of Theorem~\ref{th1} is given in
Section~\ref{sec4}. 

\section{General large deviation properties}\label{sec2}  
In this section, we recall the definition of large deviation
principle for scaled processes in $\R^d$ and in $D([0,T],\R^d)$ and some general
properties of the corresponding rate functions. 

\subsection{Large deviations}
Recall that the family scaled processes $Z^\eps(T) = \eps Z(T/\eps)$ is said to satisfy
{\em weak large deviation principle } in $\R^d$ with a rate function $I_T :\R^d\times\R^d\to\R_+$ if 
\begin{itemize}
\item[--] the function $I_T :\R^d\times\R^d\to\R_+$ is lower semicontinuous;
\item[--] for any $q\in\R^d$ and any open set $O\subset\R^d$ 
\be\label{e2-1}
\lim_{\delta\to 0}~\liminf_{\eps\to 0} ~\eps \log ~\inf_{z\in E : |\eps z - q| < \delta}
~\P_z(Z^\eps(T) \in O) ~\geq~ - \inf_{q'\in O} I_T(q,q') 
\ee
\item[--] for any $q\in\R^d$ and any compact set $V\subset\R^d$ 
\be\label{e2-2}
\lim_{\delta\to 0}~\liminf_{\eps\to 0} ~\eps \log ~\sup_{z\in E : |\eps z - q| < \delta}
~\P_z(Z^\eps(T) \in V) ~\leq~ - \inf_{q'\in V} I_T(q,q') 
\ee 
\end{itemize}
If moreover, the last inequality holds for all
closed subsets $V\subset\R^d$ then the family scaled processes $Z^\eps(T) = \eps
Z(T/\eps)$ is said to satisfy {\em large deviation principle } in $\R^d$. 

\medskip

\begin{prop}\label{pr2-1}  If the family 
$Z^\eps(T)=\eps Z(T/\eps)$ satisfies weak large deviation principle in $\R^d$ with a   
rate functions $I_{T}$  for some $T>0$, then   also it 
satisfies weak large deviation  principle in 
$\R^d$ for any
$T>0$ and the  rate function $I_T$ satisfies the following relations : 
\be\label{e2-3}
I_{\theta T}(\theta q,\theta q') ~=~ \theta ~I_{T}(q, q'), \quad \quad \forall ~\theta >0,
~T>0, ~q,q'\in\R^d 
\ee
and 
\be\label{e2-4}
I_ {T+T'}(q,q'') ~\leq~ I_T(q,q') + I_{T'}(q',q'') \quad \quad \forall ~T>0,~T'>0, ~q,q',q''\in\R^d 
\ee
\end{prop} 
\begin{proof}The first assertion of this proposition and the equality \eqref{e2-3} follows
  from contraction principle  and the identity  $Z^{\eps \theta}(t) =  \theta
  Z^{\eps}(t/\theta)$  because the mapping 
  $q \to \theta q$ is homeomorphic. Relation \eqref{e2-4} is a consequence of Markov
  property. 
\end{proof}

\subsection{Sample path large deviations}
Let $D([0,T],\R^{d})$ denote the set of all right continuous with left
limits functions from $[0,T]$ to $\R^{d}$ endowed with Skorohod metric
(see Billingsley~\cite{Billingsley}). Recall that a mapping $I_{[0,T]}:~D([0,T],\R^{d})\to
[0,+\infty]$ is called a good rate function on $D([0,T],\R^{d})$ if for any $c\geq 0$ 
and  any compact set $V\subset \R^{d}$, the set
\[
\{ \varphi \in D([0,T],\R^{d}): ~\phi(0)\in V \; \mbox{
and } \; I_{[0,T]}(\varphi) \leq c \}
\]
is compact in $D([0,T],\R^{d})$.  According to this definition, a good
rate function is lower semi-continuous. 

The family of scaled processes $(Z^\eps(t),
\,t\in[0,T])$,  is said to
satisfy {\it sample path large deviation principle} in $D([0,T], \R^{d})$ with a rate function
$I_{[0,T]}$ if for any $z\in\R^{d}$ 
\begin{equation}\label{e2-5}
\lim_{\delta\to 0} \;\liminf_{\eps\to 0} \; \inf_{z'\in E : |\eps z'-z|<\delta} \eps
\log\P_{z'}\left( Z^\eps(\cdot)\in {\cal 
O}\right) ~\geq~ -\inf_{\phi\in{\cal O}:\phi(0)=z} I_{[0,T]}(\phi)  
\end{equation}
for every open set ${\cal
O}\subset D([0,T],\R^{d})$,
and
\begin{equation}\label{e2-6}
\lim_{\delta\to 0} \;\limsup_{\eps\to 0} \; \sup_{z' \in  E : |\eps z'-z|< \delta} \eps\log\P_{z'}\left( Z^\eps(\cdot)\in
F\right) ~\leq~ -\inf_{\phi\in F:\phi(0)=z} I_{[0,T]}(\phi) 
\end{equation}
 for every closed set $F\subset   D([0,T],\R^{d})$. 

 We refer to sample path large deviation
principle as SPLD principle. Inequalities (\ref{e2-5}) and (\ref{e2-6}) are
referred as lower and upper SPLD bounds respectively.

Contraction principle applied with the continuous mapping $\phi \to \phi(T)$ from
$D([0,T],\R^d)$ to $\R^d$ proves the following statement. 

\begin{prop}\label{pr2-2} Suppose that the family of scaled processes
$(Z^\eps(t), t\in[0,T])$ satisfies SPLD principle in $D([0,T],\R^d)$ with a   
good rate functions $I_{[0,T]}$ , then the family $Z^\eps(T)$ satisfies large deviation principle in
$\R^d$ with the rate function 
\[
I_T(q,q') ~=~ \inf_{\phi(0)=q, \, \phi(T)=q'} ~I_{[0,T]}(\phi) 
\]
where the infimum is taken over all $\phi\in D([0,T],\R^d)$ with given $\phi(0)=q$ and
$\phi(T)=q'$. 
\end{prop}

\subsection{Quasipotential} The quantity 
\[
I(q,q') ~=~ \inf_{T>0} ~I_T(q,q') 
\]
represents the optimal large deviation cost to go from $q$ to $q'$. Following Freidlin and
Wentzell terminology~\cite{W-F}, such a function $I : \R^d\times\R^d\to \R_+$ is called
{\em quasipotential}. Borovkov 
and Mogulskii~\cite{Borovkov:01} called this function {\em second deviation rate function}. 

Proposition~\ref{pr2-1} implies the following properties of the function $I(q,q')$. 
\begin{cor}\label{cor2-1} If a family of scaled processes
$Z^\eps(T)=\eps Z(T/\eps)$ satisfies weak large deviation principle in $\R^d$ with a   
rate functions $I_T$ , then 
\be\label{e2-7}
I(q,q') ~=~ \inf_{T>0}  ~T ~I_1(q/T,q'/T), \quad \quad \forall~ q,q'\in\R^d,   
\ee
\be\label{e2-8}
I(\theta q,\theta q') ~=~ \theta I(q,q'), \quad \quad \forall~ \theta >0, ~q,q'\in\R^d,  
\ee
and 
\be\label{e2-9}
I(q,q') + I(q',q'') ~\geq~ I(q,q'') \quad \quad \forall~ q,q',q''\in\R^d.
\ee
\end{cor}

\section{Preliminary results}\label{sec3}

\begin{lemma}\label{lem3-1} 
Under the hypotheses (H2),  for any $q\not=q'$, $q,q'\in\R^d$ and  any $A>0$, 
   there is  $\kappa>0$, such that 
\be\label{3-1}
\lim_{\delta\to 0} ~\limsup_{n\to\infty} ~\sup_{z\in E : |z-nq|<\delta
  n} ~\frac{1}{n}~\log~\int_0^{\kappa n} \P_{z}\bigl(Z(t) \in n B(q',\delta)\bigr) \, dt ~\leq~ - A.
\ee
\end{lemma}
\begin{proof}  By Chebychev's
 inequality, for any $c>0$, any $t\in\R_+$, 
 any $a\in\R^d$ satisfying the inequality $|a|\leq c$ and any $z\in E$ satisfying the inequality
 $|z-nq|<\delta n$,  the following inequality holds~: 
\begin{align*}
 \P_{z}\bigl(|Z(t) - nq'| < \delta n\bigr) &~\leq~ \exp(- a\cdot q' n + c \delta n) ~\E_z\left(\exp(a\cdot
  Z(t))\right) \\ &~\leq~ \exp(-  a\cdot q' n + c \delta n +  a\cdot z)  ~ M_c^{t + 1}
\\ &~\leq~ \exp(-  a\cdot (q'-q) n + 2 c\delta  n)  ~M_c^{t +1}
\end{align*}
where the constant 
\[
M_c ~\dot=~ \sup_{a\in\R^d : |a| \leq c} ~\sup_{z\in E} ~\sup_{t\in[0,1]} \E_z\left(e^{ a\cdot
  (Z(t)-z)}\right)  ~=~ \sup_{a\in\R^d : |a| \leq c}  \hat\varphi(a)  ~\geq~ 1
\]
is finite because of the Assumption (H2). 
Hence, letting $c ~=~ 2 A/|q'-q|$, $a=2 A (q'-q)/|q'-q|^2$ and $\kappa =
 A / (2\log M_c)$, for any $z\in E$ satisfying the inequality $|z-nq|<\delta n$,   
we obtain 
\begin{align*}
\int_0^{\kappa n}  \P_{z}\bigl(|Z(t) - nq'|<\delta n\bigr) \, dt 
&\leq~ \exp\bigl(- c |q'-q| n + (\delta +\eps) c n \bigl) \int_0^{\kappa n} 
M^{t+1}_c \, dt \\ 
&\leq~ \exp\left(- 2 A n + c (\delta+\eps) n + (\kappa n +1) \log M_c\right)/\log M_c \\
&\leq~ \exp(-  A n) M_c/\log M_c 
\end{align*}
whenever $0 < \delta  < A/(4c)$. The last inequality proves \eqref{3-1}. 
\end{proof}

\begin{lemma}\label{lem3-2} Suppose that   the upper large deviation bound \eqref{e2-2}
  holds for every $q\in{\cal R}$ and for every compact set $V\subset \R^d$ with a rate
  function $I_T: \R^d\times\R^d\to \R_+$ satisfying the inequality \eqref{e2-3}. 
Then for any $0 < \kappa < K < \infty$
  and any $q,q'\in\R^d$ 
\be\label{3-2}
\lim_{\delta\to 0} ~\limsup_{n\to\infty} ~\sup_{z\in E : |z-nq|<
  \delta n} ~\frac{1}{n}~\log~\int_{\kappa n}^{K n} \P_{z}\bigl(|Z(t) - nq'| < \delta n 
  \bigr) \, dt ~\leq~ - I(q,q').
\ee
\end{lemma}
\begin{proof} 
The upper large deviation bound \eqref{e2-2} applied with $\eps=1/t$ proves that for any $x,y\in\R^d$, 
\[
\lim_{\delta\to 0}\limsup_{t\to\infty} ~\sup_{z : |z - t x|<\delta t} ~\frac{1}{t} \log \P_z\left( |Z(t)/t - y| <\delta 
\right) ~\leq~  - I_1(x,y) ~\leq~ - I(x,y).
\]
Letting 
\[
I^\sigma(q,q') ~\dot=~ \min\{I(q,q'), (|q|+|q'|)/\sigma\}
\]
for $\sigma >0$ and using the above inequality with $x=\theta q$ and $y=\theta q'$ we
get that  for any 
$\sigma>0$ and $\theta >0$, there are $0<\delta(\theta)<\sigma/K$ and 
$t(\theta)>0$ such that 
for all $t > t(\theta)$ and  $0<\delta\leq\delta(\theta)$, the following inequality holds 
\begin{align*}
\sup_{z : |z - t\theta q|<\delta t} \P_z\left( |Z(t) - t \theta q'| <\delta t
\right) &~\leq~ \exp(- t I^\sigma(\theta q,\theta q') + \sigma t) \\ &~\leq~  \exp(- t
\theta I^\sigma(q,q') +\sigma t) 
\end{align*}
(the last relation is a consequence of \eqref{e2-3}). Since the set $[K^{-1
  },\kappa^{-1}]$ is compact then there are 
$\theta_1,\ldots,\theta_m\in [K^{-1},\kappa^{-1}]$ such that 
\[
[K^{-1},\kappa^{-1}] \subset \bigcup_{i=1}^m ]\theta_i - \delta(\theta_i)/c, \theta_i +
  \delta(\theta_i)/c[ 
\]
with $c=2\max\{1,|q|,|q'|\}$. Denote $t_i = t(\theta_i)$
and $\delta_i = \delta(\theta_i)$, and let  
\[
\delta = \kappa \min_{1\leq i\leq m}  
\delta_i/2.  
\]
Then for every 
$t\in[\kappa n, K n]$ there is $i(t)\in\{1,,\ldots,m\}$ such that 
\be\label{3-3}
|n  -
t \theta_{i(t)}|<\delta_{i(t)} t /2,  
\ee
\[|n q  -
\theta_{i(t)} t q |< \delta_{i(t)}t/2 \quad \text{ and } \quad |n q' -
\theta_{i(t)} t q' |<  \delta_{i(t)}t/2 
\]
For any $n\geq 
\kappa^{-1}\max\{t_1,\ldots,t_m\}$,  $t\in[\kappa n, K n]$ and  $z\in E$ satisfying
the inequality $|z-n q| < \delta n$, we obtain therefore  
\begin{align}
 \P_{z}\left(|Z(t)- n q'  | < \delta n\right) \nonumber &~\leq~ \sup_{z : |z-n q|<\delta
  t/\kappa} \P_{z}\left(\left|Z(t) - n q' \right| < \delta t/\kappa\right) \nonumber\\ 
  &~\leq~ \sup_{z : |z- \theta_{i(t)} t q|<\delta_{i(t)} t}  \P_{z}\left(|Z(t) - 
  \theta_{i(t)} t q'|
  < \delta_{i(t)} t\right) \nonumber\\  &~\leq~ \exp(- t\theta_{i(t)} I^\sigma(q,q') +
  \sigma t) \label{3-4}
\end{align}
Moreover,  \eqref{3-3} shows that  for every $t\in[\kappa n, K n]$,  
\begin{align*}
 t\theta_{i(t)} ~\geq~ (n - t\delta_{i(t)}) ~\geq~ n (1 - K\delta_{i(t)}) ~\geq~ n (1 - \sigma) 
\end{align*}
The last relation combined with \eqref{3-4} proves
that for any $n\geq 
\kappa^{-1}\max\{t_1,\ldots,t_m\}$, 
\[
\int_{\kappa n}^{Ka} \P_{z}\left(|Z(t)- n q'  | < \delta n\right)  \, dt ~\leq~ (K-\kappa)n ~\exp\left(-
n (1 - \sigma) I^\sigma(q,q')
+ K n\sigma\right)  
\] 
and consequently,
\[
~\limsup_{n\to\infty} \sup_{z\in E : |z-nq|<
  \delta n} ~\frac{1}{n}~\log~\int_{\kappa n}^{K n} \! \P_{z}\bigl(|Z(t) - nq'| < \delta n 
  \bigr) \, dt ~\leq~ - (1 - \sigma) I^\sigma(q,q') + K \sigma .
\]
Letting finally $\delta \to 0$ and $\sigma\to 0$ we get \eqref{3-2}. 
\end{proof}

\begin{lemma}\label{lem3-3} If  the upper large deviation bound \eqref{e2-2} holds with a rate
  function $I_T$ such that $I_T(0,0) >0$ then for any 
  $A>0$, any $q\in\R^d$ and any bounded set $V\subset\R^d$, there is  
  $K>0$ such that 
\be\label{3-5}
\lim_{\delta\to 0} ~\limsup_{n\to\infty} ~\sup_{z\in E : |z-nq|<
  \delta n} ~\frac{1}{n}~\log~\int_{K n}^\infty \P_{z}\bigl(Z(t) \in nV\bigr) \, dt ~\leq~ - A.
\ee
\end{lemma}
\begin{proof} From \eqref{e2-2} it follows that 
\begin{align*}
\lim_{\delta\to 0}\limsup_{t\to\infty} ~\sup_{z : |z| < \delta t} ~\frac{1}{t} \log &\P_z\left( |Z(t)| \leq \delta
t\right) \\ &~=~ \lim_{\delta\to 0}\limsup_{\eps\to 0} ~\sup_{z : \eps |z| < \delta T}
~\frac{\eps}{T} \log \P_z\left( |Z^\eps(T)| \leq \delta 
T\right) \\ &~\leq~ - \frac{1}{T} ~\lim_{\delta\to 0} ~\inf_{z' : |z'|\leq\delta}
~I_T(0,z') ~=~ - \frac{1}{T} I_T(0,0). 
\end{align*}
Hence, there are $\delta_0 > 0$ and $t_0>0$ such
  that for  all $t>t_0$ 
\[
~\sup_{z : |z| < \delta_0 t} \P_z\left( |Z(t)| \leq \delta_0
t\right) ~\leq~ \exp\left(- I_T(0,0) t/(2T)\right).
\]
For  $t> \max\{t_0, (|q| + 1)a/\delta_0, \sup_{q'\in V} |q'| a/\delta_0\}$ we get
therefore 
\begin{align*}
\sup_{z : |z-n q|< a} \P_{z}\bigl(Z(t)\in n V\bigr)  &~\leq~ ~\sup_{z :
  |z| < \delta_0 t} \P_z\left( |Z(t)| \leq \delta_0
t\right) \\ &~\leq~ \exp\left(- I_T(0,0) t/(2T)\right).
\end{align*}
This inequality shows that for  any 
$K> \max\{(|q| + 1)/\delta_0, \sup_{q'\in V} |q'| /\delta_0\}$ and  $n > t_0/K$,  
\begin{align*}
\sup_{z : |z-n q|< n} \int_{Kn}^\infty \P_{z}\bigl(Z(t)\in n V\bigr) \, dt &~\leq~
~\int_{K n}^\infty \exp\left(- I_T(0,0) t/(2T)\right) \, dt \\ &~\leq~ \frac{ 2 T \exp\left(-
  I_T(0,0) Kn/(2T)\right)}{I_T(0,0)} 
\end{align*}
Letting $K = \max\{(|q|+ \sup_{q'\in V} |q'|+ 1)/\delta_0, ~2 A T/ I_T(0,0)\}$, we get 
\[
\sup_{z : |z-n q|<\delta n} \int_{Kn}^\infty  \P_{z}\bigl(Z(t)\in n V\bigr) \,
dt ~\leq~ 2T \exp(-A n)/ I_T(0,0) 
\]
for all $n >t_0/K$  and  $0<\delta < 1$, 
and consequently, \eqref{3-5} holds.  
\end{proof}

\section{Proof of Theorem~\ref{th1}}\label{sec4} 
\subsection {Proof of the assertion (i)} For discrete time Markov process $(Z(t))$, the
proof of this assertion is straightforward : if the
  family of scaled processes $Z^\eps (T) = \eps Z(T/\eps)$ satisfy the lower large deviation
  bound  with a rate function $I_T$ then  for any $q\in{\cal R}$ and any open set $O\subset\R^d$, 
\begin{align*}
\lim_{\delta\to 0} ~\liminf_{n\to \infty}  \inf_{z\in E : |z-nq| < 
  \delta n} &~\frac{1}{n}\log  ~G(z, n O) \\ &~\geq~ 
\lim_{\delta\to 0} ~\liminf_{n\to \infty} ~\inf_{z\in E : |z-nq| < 
  \delta n} ~\frac{1}{n}\log P_z(Z(nT) \in n O) \\
&~\geq~ 
\lim_{\delta\to 0} ~\liminf_{\eps\to 0} ~\inf_{z\in E : |\eps z-q| < 
  \delta } ~\eps \log P_z(Z^\eps(T) \in O) \\
&~\geq~  - \inf_{q'\in O} I_T(q,q').  
\end{align*}
The last inequality proves  the lower bound \eqref{e1-1} because  $T>0$ is arbitrary. \\

Suppose now that $(Z(t))$ is a continuous time Markov process and let us show that for any
$q,q'\in \R^d$ and $\delta,\delta' >0$, the following inequality holds 
\begin{align}
\liminf_{n\to \infty}  \inf_{z\in E : |z-nq| < 
  \delta n} &~\frac{1}{n}\log G\bigl(z, n B(q',\delta')\bigr) \nonumber\\ &~\geq~
  \liminf_{n\to \infty}  \inf_{z\in E : |z-nq| <   
  \delta n} ~\frac{1}{n}\log ~\P_z\bigl(Z(Tn)\in n B(q',\delta'/2)\bigr).  \label{e4-1p}
\end{align}
Indeed, for any $\delta > 0$,
$a\in\R^d$, $t\in [0,1]$ and $z\in E$, by Chebyshev's inequality 
\[
\P_z \bigl(a\cdot (Z(t) - z) \geq n\delta\bigr) ~\leq~ e^{-\delta n} \E_z\left(e^{a\cdot (Z(t)
  -z)}\right) ~\leq~  e^{-\delta n} ~\hat\varphi(a) 
\]
from which it follows that 
\begin{align*}
\P_z \Bigl(|Z(t) - z|\geq n\delta\Bigr) \leq~ \sum_{a\in\Z^d : |a|=1} \P_z \left(a\cdot
(Z(t) - z) \geq \frac{n\delta}{2d}\right) \leq~ 2d M_1 ~e^{-\delta n/(2d)} 
\end{align*}
where 
\[
M_1 ~\dot=~  \sup_{a :
  |a|=1} ~\hat\varphi(a).   
\]
The last inequality shows that 
\begin{align*}
\sup_{z\in E : |z- q' n|  < \delta n/2} ~\P_z\bigl(|Z(t) - q'| \geq \delta n\bigr)  &~\leq~
\sup_{z\in E : |z- q' n|  < \delta n/2} ~\P_z\bigl(|Z(t) - z| \geq \delta n/2\bigr) \\
&~\leq~ 2d M_1 ~e^{-\delta n/(4d)} 
\end{align*}
or equivalently, that 
\[
\inf_{z\in E : |z- q' n|  < \delta n/2} ~\P_z\bigl(|Z(t) - q'| < \delta n\bigr) ~\geq~ 1 -
2d M_1 ~e^{-\delta n/(4d)}. 
\]
Hence, for $t\in[Tn, Tn +1]$,  using Markov property we obtain    
\begin{align*}
\P_z\bigl(Z(t)\in n B(q',\delta')\bigr) &~\geq~
\P_z\bigl(Z(Tn)\in n B(q',\delta'/2)\bigr)  \\ &\quad \quad \quad \quad \quad \quad \times 
\inf_{z\in E : |z- q' n|  < \delta n/2} ~\P_z\bigl(|Z(t-Tn) - q'| < \delta' n\bigr)
\\&~\geq~  \P_z\bigl(Z(Tn)\in n B(q',\delta'/2)\bigr)  \left(1 - 2d  M_1 ~e^{-\delta' n/(4d)} \right)
\end{align*}
and consequently, 
\begin{align*}
G\bigl(z, n B(q',\delta')\bigr) &~\geq~ \int_{Tn}^{Tn+1} \P_z\bigl(Z(t)\in n
B(q',\delta')\bigr)  \, dt \\ &~\geq~  \P_z\bigl(Z(Tn)\in n B(q',\delta'/2)\bigr)  \left(1
- 2d  M_1 ~e^{-\delta' n/(4d)} \right) 
\end{align*} 
The above inequality proves \eqref{e4-1p} because 
\[
\lim_{n\to\infty} \frac{1}{n} ~\log  \left(1
- 2d  M_1 ~e^{-\delta' n/(4d)} \right)  = 0. 
\]
Now, using the inequality \eqref{e4-1p} and the large deviation lower bound \eqref{e-lb}
with an open set $B(q',\delta'/2)$  
we conclude that for any open set $O\subset \R^d$
and for any $q'\in O$ and $\delta' > 0$ such that $B(q',\delta') \subset O$, 
\begin{align*}
\lim_{\delta\to 0} ~\liminf_{n\to \infty}  &\inf_{z\in E : |z-nq| < 
  \delta n} ~\frac{1}{n}\log  ~G\bigl(z, n O\bigr) \\&~\geq~ \lim_{\delta\to 0}
  ~\liminf_{n\to \infty}  \inf_{z\in E : |z-nq| <  
  \delta n} ~\frac{1}{n}\log  ~G\bigl(z, n B(q',\delta')\bigr) \\ 
& ~\geq~ \lim_{\delta\to 0} ~\liminf_{n\to \infty}  \inf_{z\in E : |z-nq| < 
  \delta n} ~\frac{1}{n}\log  ~\P_z\bigl(Z(Tn) \in  nB(q',\delta'/2)\bigr) \\ 
& ~\geq~ \lim_{\delta\to 0} ~\liminf_{n\to \infty}  \inf_{z\in E : |z-nq| < 
  \delta n} ~\frac{1}{n}\log  ~\P_z\bigl(Z^{1/n}(T) \in B(q',\delta'/2)\bigr) \\ &~\geq~ -
  \inf_{q''\in B(q',\delta'/2)} I_T(q,q'') ~\geq~ - 
  I_T(q,q') 
\end{align*}
The last inequality proves the lower bound \eqref{e1-1} because $T>0$ and $q'\in O$ are arbitrary. 

\subsection{Proof of the assertion (ii)} To prove the upper bound \eqref{e1-2} it is
sufficient to show that for any $q,q'\in{\cal R}$, 
\be\label{e4-1}
\lim_{\delta\to 0} ~\limsup_{n\to\infty} ~\sup_{z\in E : |z-nq|<\delta 
  n} ~\frac{1}{n}\log  G\bigl(z, n B(q',\delta)\bigr) ~\leq~ -  \hat{I}(q,q') 
\ee
With such a local upper large deviation bound one can obtain the upper bound \eqref{e1-2}
by using exactly 
the same arguments as in  the proof of Theorem~4.1.11 in the book of Dembo and
Zeitouni~\cite{D-Z}.  

 For $q=q'$, the local upper bound \eqref{e4-1} is a
consequence of Lemma~\ref{lem3-2} applied with $V= B(q',\delta)$ because  
\begin{align*}
\lim_{\delta\to 0} ~\limsup_{n\to\infty} ~\sup_{z\in E : |z-nq|<\delta
  n} ~\frac{1}{n}\log~\int_0^{Kn} &\P_{z}\bigl(|Z(t) - 
  nq'|<\delta n\bigr) \, dt \\ &\leq~ \lim_{\delta\to 0} ~\limsup_{n\to\infty} \frac{1}{n}~\log~ Kn
  ~\leq~ 0 ~=~ \hat{I}(q,q).
\end{align*}
\medskip  
Suppose now that $q\not=q'$ and let us choose 
 the constants $K>\kappa>0$ for which the inequalities \eqref{3-1} and \eqref{3-5} hold 
with $V = B(q',\delta)$ and 
\[
A ~>~ I^\sigma(q,q') ~\dot=~ \min\{I(q,q'), (|q|+|q'|)/\sigma\}
\]
for some $\sigma>0$. Then using these inequalities together with Lemma~\ref{lem3-2} and
Lemma~1.2.15 of~\cite{D-Z} we get 
\[
\lim_{\delta\to 0} ~\limsup_{n\to\infty} ~\sup_{z\in E : |z-nq|<\delta 
  n} ~\frac{1}{n}\log  G\bigl(z, n B(q',\delta)\bigr) ~\leq~ -   I^\sigma(q,q') 
\]
and hence, letting $\sigma\to 0$ we obtain \eqref{e4-1}.

\subsection{Proof of the assertion (iii)} 

To prove the inequality \eqref{e1-3} we use  SPLD upper
bound, Lemma~\ref{lem3-3} and the inequality 
\begin{align*}
G(z, n V) - G_{nR}(z, n V) &~\leq~ \int_0^{Tn} \P_z\bigl(\tau_{nR} \leq t\bigr) \, dt ~+~ \int_{T
  n}^\infty \P_{z}\bigl(Z(t) \in n V\bigr) \, dt \\
&~\leq~  T n ~\P_z\bigl(\tau_{nR} \leq T n\bigr) ~+~ \int_{T
  n}^\infty \P_{z}\bigl(Z(t) \in n V\bigr) \, dt. 
\end{align*}
If the family of scaled processes $(Z^\eps(t), t\in
[0,T])$ satisfies SPLD upper bound with a good rate function $I_{[0,T]}$ then for any
$A>0$, the set $$\{\phi \in 
D([0,T],\R^d) : \phi(0) = q, \, I_{[0,T]}(\phi)\leq A\}$$ is compact and therefore
bounded. Hence, 
  there is $R>0$ for which 
  the infimum of $I_{[0,T]}(\phi)$ over all $\phi \in D([0,T],\R^d)$ with $\phi(0)=q$ and
  $\sup_{s\in[0,T]}|\phi(s)| \geq R$ is greater than $A$. Using
  SPLD upper bound with $\eps=1/n$ from this it follows that   
\[
\lim_{\delta\to 0}\limsup_{n\to\infty} ~\sup_{z : |nq - z|<\delta n} ~\frac{1}{n} \log
\P_{z}\left( \sup_{s\in[0,Tn]} |Z(s)| >  R n
\right) ~\leq~  - A. 
\]
Hence, under the hypotheses (H1'), for any $T>0$ and any $A>0$ there exists $R>0$ such
that 
\[
\lim_{\delta\to 0}\limsup_{n\to\infty} ~\sup_{z : |n q - z|<\delta n} ~\frac{1}{n} \log
\P_{z}\bigl( \tau_{nR} \leq T n \bigr) ~\leq~ - A 
\]
Choosing $T>0$ in such a way that for a given $A>0$, 
\[
\lim_{\delta\to 0}\limsup_{n\to\infty} ~\sup_{z : |n q - z|<\delta n} ~\frac{1}{n} \log
 \int_{Tn}^\infty \P_{z}\bigl(Z(t) \in n V\bigr) \, dt ~\leq~ - A 
\] 
(the existence of such $T>0$ is proved by Lemma~\ref{lem3-3}) and using Lemma~1.2.15 of
\cite{D-Z} we get \eqref{e1-3}. 

\subsection{Proof of the assertion (iv)} 
If the Markov
process $(Z(t))$ satisfies the weak communication condition (H3), then  clearly 
\[
\log G(z,B(z',\sigma(z'))) ~\geq~ - \theta |z-z'| 
\]
for all $z,z'\in E$. For $q\not=q'$, $q,q'\in {\cal R}$, using 
the local upper bound \eqref{e4-1} from this it follows that  
\[
- \theta |q-q'| ~\leq~ \lim_{\delta \to 0} \limsup_n ~\sup_{z : |nq - z|<\delta n} ~\frac{1}{n} ~\log G(z,
  B(n q',\delta n)) ~\leq~ - I(q,q')  
\]
and consequently, 
\[
I(q,q') ~\leq~ \theta |q-q'|, \quad \quad \forall q\not= q', \; q,q'\in{\cal R}. 
\]
The last inequality combined with \eqref{e2-9} shows that 
\[
I(q,q) ~\leq~ I(q,q') + I(q',q) ~\leq~ 2\theta |q-q'|, \quad \quad \forall \; q\not= q', \; q,q'\in{\cal R}.
\]
Letting therefore $q'\to q$ we conclude that $I(q,q)=0$ for every $q\in{\cal R}$. 
Moreover, for any $q,q',w,w'\in{\cal R}$ we get 
\[
I(q,w) \leq I(q,q') + I(q',w') + I(w',w) \leq I(q',w') + \theta |q-q'| + \theta |w-w'|
\]
and 
\[
I(q',w') \leq I(q',q) + I(q,w) + I(w,w') \leq I(q,w) + \theta |q-q'| + \theta |w-w'|.
\]
These inequalities show that  the function $(q,q') \to I(q,q')$ is continuous on ${\cal 
  R}\times{\cal R}$.  Theorem~\ref{th1} is therefore proved.

\bibliographystyle{amsplain}
\bibliography{ref}
\end{document}